\documentclass[11pt,a4paper]{article}

\usepackage[T1]{fontenc}
\usepackage[margin=25mm,headheight=15pt]{geometry}
\usepackage{amsmath,amsthm,mathtools}
\usepackage{newtxtext,newtxmath}
\usepackage{microtype,booktabs,listings,textcomp}
\usepackage{fancyhdr}
\usepackage[table]{xcolor}
\usepackage[hidelinks]{hyperref}

\numberwithin{equation}{section}
\newtheorem{theorem}{Theorem}[section]
\newtheorem{proposition}[theorem]{Proposition}
\newtheorem{lemma}[theorem]{Lemma}
\theoremstyle{remark}

\newcommand{\Z}{\mathbb Z}

\newcommand{\lammax}{\lambda_{\max}}
\newcommand{\norm}[1]{\left\lVert#1\right\rVert}
\newcommand{\ip}[2]{\left\langle#1,#2\right\rangle}
\newcommand{\floor}[1]{\left\lfloor#1\right\rfloor}
\newcommand{\ceil}[1]{\left\lceil#1\right\rceil}
\newcommand{\eps}{\varepsilon}

\definecolor{paperblue}{RGB}{225,235,245}

\title{A counterexample to the Harper-norm silver-mean conjecture}
\author{
Evgeny Shchukin\thanks{
Johannes-Gutenberg University of Mainz, Institute of Physics, Mainz, Germany
Email: shchukin@uni-mainz.de
}
\and
Peter van Loock\thanks{
Johannes-Gutenberg University of Mainz, Institute of Physics, Mainz, Germany
Email: loock@uni-mainz.de
}
}
\date{}

\begin{document}
\maketitle
\thispagestyle{fancy}
\vspace{-2.5em}
\begin{abstract}
The silver-mean conjecture formulated in Ref.~\cite{GerholdShalit2022} states
that the norm $E(\alpha)$ of the Harper operator attains its minimum at the
silver mean $\alpha = \sqrt{2}-1$. We show that this conjecture is wrong by
constructing an explicit rational point with a smaller Harper-norm value:
\begin{displaymath}
    E(181/437)<2.591046794<E(\sqrt2-1).
\end{displaymath}
The separating decimal is an exact rational number. All acceptance checks use
exact integer arithmetic with rigorously proved coefficient enclosures. No
floating-point eigenvalue estimates or external witness data enter the
certificate.
\end{abstract}

\section{Introduction}
\label{sec:applications}

The Harper operator originates in the tight-binding description of
electrons on a two-dimensional lattice in a uniform magnetic field
\cite{Harper1955}. Its spectral dependence on the magnetic flux gives
rise to the Hofstadter butterfly \cite{Hofstadter1976}. For irrational
frequencies, it is also the critical almost Mathieu operator, a basic
example of a quasiperiodic Schr\"odinger operator. Here we study a
specific question about its norm rather than its full spectrum.

For $\alpha,\phi\in\mathbb R$, define the Harper operator on
$\ell^2(\mathbb Z)$ by
\begin{equation}\label{eq:harper}
    (H_{\alpha,\phi}x)_n
    =
    x_{n-1}+x_{n+1}
    +2\cos(2\pi\alpha n+\phi)x_n,
    \qquad n\in\mathbb Z.
\end{equation}
The frequency $\alpha$ is measured in turns and the phase $\phi$ in
radians. The operator is bounded and self-adjoint, with
$\|H_{\alpha,\phi}\|\le 4$. We define the \emph{Harper norm} by
\begin{equation}\label{eq:intro-E}
    E(\alpha)
    =
    \sup_{\phi\in\mathbb R}\|H_{\alpha,\phi}\|.
\end{equation}
For a bounded self-adjoint operator $A$, write
\begin{equation}
    e(A)
    =
    \sup\operatorname{spec}(A)
    =
    \sup_{x\ne0}
    \frac{\langle x,Ax\rangle}{\|x\|^2}
\end{equation}
for its upper spectral edge. The unitary operator
$(Jx)_n=(-1)^n x_n$ satisfies
\begin{equation}
    JH_{\alpha,\phi}J=-H_{\alpha,\phi+\pi}.
\end{equation}
Consequently, the union of the phase spectra is symmetric about
zero, and
\begin{equation}\label{eq:intro-phase-edge}
    E(\alpha)
    =
    \sup_{\phi\in\mathbb R}e(H_{\alpha,\phi}).
\end{equation}
The quantity $E(\alpha)$ agrees with the universal rotation-algebra
norm considered by Gerhold and Shalit
\cite{GerholdShalit2022}. In their notation,
\begin{equation}
    E(\alpha)=\|h_{2\pi\alpha}\|,
    \qquad
    h_\theta=u_\theta+u_\theta^*+v_\theta+v_\theta^*,
\end{equation}
where $u_\theta,v_\theta$ are the universal unitaries satisfying
$v_\theta u_\theta=e^{i\theta}u_\theta v_\theta$.
They proved that the optimal constant $c_\theta$ for dilating
$e^{i\theta}$-commuting unitaries to scalar multiples of commuting
unitaries satisfies
\begin{equation}\label{eq:intro-dilation}
    c_{2\pi\alpha}=\frac{4}{E(\alpha)}.
\end{equation}
Their numerical experiments led them to conjecture that $c_\theta$ attains its
global maximum at $\theta=2\pi(\sqrt2-1)$. Equivalently, writing $s_*=\sqrt2-1$,
the \emph{silver-mean conjecture} asserts that
\begin{equation}\label{eq:intro-conjecture}
    E(s_*)\leqslant E(\alpha)
    \qquad\text{for every }\alpha\in\mathbb R.
\end{equation}
The periodicity $E(\alpha+1)=E(\alpha)$ and symmetry $E(-\alpha)=E(\alpha)$
allow the minimization to be restricted to $\alpha\in[0,1/2]$. We disprove this
conjecture by establishing
\begin{equation}\label{eq:intro-counterexample}
    E(181/437)
    <
    2.591046794
    <
    E(s_*).
\end{equation}
The separating decimal denotes the exact rational number
$25910467940/10^{10}$. The candidate is discovered numerically,
but the proof of \eqref{eq:intro-counterexample} does not use
floating-point eigenvalue estimates.

For the finite approximations used below, it is useful to record
the stronger zero-phase identity
\begin{equation}\label{eq:intro-zero-phase}
    E(\alpha)
    =
    e(H_{\alpha,0})
    =
    \|H_{\alpha,0}\|,
    \qquad \alpha\in\mathbb R.
\end{equation}
For irrational $\alpha$, the function
$\phi\mapsto e(H_{\alpha,\phi})$ is continuous, since
\begin{equation}
    |e(H_{\alpha,\phi})-e(H_{\alpha,\psi})|
    \leqslant
    \|H_{\alpha,\phi}-H_{\alpha,\psi}\|
    \leqslant 2|\phi-\psi|.
\end{equation}
Translation of the sequence index makes this function invariant under phase
shifts by $2\pi\alpha$. These shifts have dense orbits modulo $2\pi$, so the
function is constant. Equation~\eqref{eq:intro-phase-edge} then gives
$e(H_{\alpha,0})=E(\alpha)$.

For rational $\alpha=p/q$, Gerhold--Shalit's finite-dimensional reduction
identifies the universal norm with the largest eigenvalue $\lambda$ of the
zero-phase periodic $q\times q$ matrix \cite{GerholdShalit2022}; this matrix is
denoted $H_{p/q}$ in Section~2. Every $H_{p/q,\phi}$ represents the universal
Harper element, so $E(p/q)\leqslant\lambda$. Conversely, repeating a unit top
eigenvector of the periodic matrix over $N$ consecutive periods and setting it
equal to zero elsewhere gives Rayleigh quotients for $H_{p/q,0}$ tending to
$\lambda$: the squared norm is $N$, whereas the boundary correction is
independent of $N$. Hence $e(H_{p/q,0})\geqslant\lambda$, proving
$e(H_{p/q,0})=E(p/q)$. In either case,
\begin{equation}
    e(H_{\alpha,0})
    \leqslant \|H_{\alpha,0}\|
    \leqslant E(\alpha),
\end{equation}
which completes the justification of
\eqref{eq:intro-zero-phase}.

The certificate consists of two finite positive sequences.
The first is closed into a periodic vector of length $437$ and
provides an upper bound for the largest eigenvalue of the rational
Harper matrix. The second is extended evenly to the interval
$[-283,283]$ and then by zero outside it; its Rayleigh quotient
provides a lower bound for $E(s_*)$. Both sequences are generated
by second-order recurrences with directed integer rounding.
Explicit rational Taylor bounds enclose the cosine coefficients,
and finite sign checks verify the required inequalities.
No approximate eigenvectors or precomputed cosine arrays are
inputs to the certificate.

Section~2 describes the numerical search. Sections~3--5 develop
the positive-vector bounds and the general certificate for
arbitrary sequence lengths. Sections~6--7 give the concrete
parameters, prove the coefficient enclosures, and present the
exact arithmetic establishing
\eqref{eq:intro-counterexample}.

\section{A candidate counterexample}

For a reduced rational $\alpha\equiv r = p/q$ define the matrix $H_{r}$ on indices
$0,\ldots,q-1$ by
\begin{equation}\label{eq:A}
 (H_{r}x)_j=x_{j-1}+x_{j+1}+2\cos(2\pi rj)x_j,
\end{equation}
where all indices are interpreted modulo $q$. For $q\geqslant 3$ it explicitly
reads as
\begin{equation}
    H_{r}
    =
    \begin{pmatrix}
    c_0(r) & 1 & 0 & \cdots & 1\\
    1 & c_1(r) & 1 & \ddots & \vdots\\
    0 & 1 & c_2(r) & \ddots & 0\\
    \vdots & \ddots & \ddots & \ddots & 1\\
    1 & \cdots & 0 & 1 & c_{q-1}(r)
    \end{pmatrix},
\end{equation}
where $c_j(\alpha)=2\cos(2\pi\alpha j)$. For $q=1$ and $q=2$ this matrix is 
\begin{equation}
    H_{p/1} = (4), \quad
    H_{p/2} = 
    \begin{pmatrix}
        2 & 2 \\
        2 & -2
    \end{pmatrix}, \quad
    (p\ \text{odd}).
\end{equation}
According to Ref.~\cite{GerholdShalit2022} we have the relation
\begin{equation}\label{eq:rational}
 E(r)=\lammax(H_{r}).
\end{equation}
It follows that for a rational $\alpha$ the Harper norm $E(\alpha)$ is an
algebraic number so, in principle, one could compute it with an arbitrary
precision with the help of a computer algebra system to get the desired
inequality (and hope that there are no bugs in the software used for this task).
We use this approach only to search for a candidate rational $r$ that is likely
to violate the conjecture. 

For an irrational $\alpha$ there is no simple way to compute $E(\alpha)$, so we
need another approach. For any real $\alpha$ we construct a sequence of numbers
that tends to $E(\alpha)$ from below. We search for a candidate rational $r$ for
which the Harper norm $E(r)$ is smaller than an element of the approximating
sequence for $E(s_*)$. When such a candidate is found, the strict inequality
$E(r)<E(s_*)$ is not automatically established since one cannot guarantee that
unavoidable roundoff errors do not destroy it. To really certify this inequality
we will use the integer arithmetic only. 

To construct an approximating sequence for $E(\alpha)$, let $B_{\alpha,m}$
denote the restriction of $H_{\alpha,0}$ to the interval $\{-m,-m+1,\ldots,m\}$
with Dirichlet boundary conditions, $m\geqslant 0$. Thus $B_{\alpha,m}$ is the
$(2m+1)\times(2m+1)$ real symmetric matrix whose quadratic form is
\begin{equation}
    x^T B_{\alpha,m}x
    =
    \sum_{j=-m}^{m}
        2\cos(2\pi\alpha j)x_j^2
    +
    2\sum_{j=-m}^{m-1}x_jx_{j+1}.
\end{equation}
Equivalently, vectors on $[-m,m]$ are extended by zero outside this interval
before applying $H_{\alpha,0}$. Since $c_{-j}(\alpha) = c_j(\alpha)$, the matrix
$B_{\alpha,m}$ commutes with the reflection $(Rx)_j=x_{-j}$. Its largest
eigenvalue has an even eigenvector. Indeed, the off-diagonal entries of
$B_{\alpha,m}$ are nonnegative, so replacing a real vector $x$ by $(|x_j|)$
cannot decrease its Rayleigh quotient. A top eigenvector may therefore be chosen
nonnegative. Averaging it with its reflection gives a nonzero even top
eigenvector.

Let $e_j$ denote the coordinate vector equal to one at site $j$ and zero at all
other sites. The even subspace has the orthonormal basis
\begin{equation}
    f_0 = e_0,\qquad
    f_j=\frac{e_j+e_{-j}}{\sqrt2},
    \qquad 1\leqslant j\leqslant m.
\end{equation}
In this basis the restriction of $B_{\alpha,m}$ to the even subspace
is the $(m+1)\times(m+1)$ tridiagonal matrix
\begin{equation}
    C_{\alpha,m}
    =
    \begin{pmatrix}
    c_0(\alpha) & \sqrt2 & 0 & \cdots & 0\\
    \sqrt2 & c_1(\alpha) & 1 & \ddots & \vdots\\
    0 & 1 & c_2(\alpha) & \ddots & 0\\
    \vdots & \ddots & \ddots & \ddots & 1\\
    0 & \cdots & 0 & 1 & c_m(\alpha)
    \end{pmatrix},
\end{equation}
For $m=0$ we put $C_{\alpha,0}=[2]$. For $m\geqslant1$, the factor $\sqrt2$ in
the first off-diagonal entry comes from the two bonds $0\leftrightarrow1$ and
$0\leftrightarrow-1$:
\begin{equation}
    \langle
        f_0,
        B_{\alpha,m}f_1
    \rangle = 
    \left\langle
        e_0,
        B_{\alpha,m}
        \frac{e_1+e_{-1}}{\sqrt2}
    \right\rangle
    =\sqrt2.
\end{equation}
All subsequent off-diagonal entries are equal to $1$.

\begin{proposition}\label{prop:dirichlet-convergence}
For every $\alpha\in\mathbb R$ and every $m\geqslant0$,
\begin{equation}
    \lambda_{\max}(C_{\alpha,m})
    =
    \lambda_{\max}(B_{\alpha,m})
    \leqslant E(\alpha).
\end{equation}
Moreover
\begin{equation}
    \lambda_{\max}(C_{\alpha,m})
    \nearrow E(\alpha)
    \qquad\text{as }m\to\infty.
\end{equation}
\end{proposition}
\noindent The proof of this proposition is given in the Appendix.

Now we can search for a possible candidate that violates the silver mean
conjecture. For a reduced rational $r = p/q$, the Harper norm $E(r)$ is the
largest eigenvalue of the Harper matrix $H_r$, which makes it inexpensive to
scan many rational frequencies. We enumerated all reduced fractions
\begin{equation}
   0<r=\frac pq<\frac12,
\end{equation}
with $3\leqslant q \leqslant 500$ in increasing order of $q$ and then $p$.  For
each fraction we computed the largest eigenvalue of $H_{r}$. During the
enumeration we kept track of successive record minima. The resulting list is
shown in the left part of Table~\ref{tbl:candidate}. The right half of that
table shows the values of $\lambda_{\max}(C_{s_*,m})$ for an increasing sequence
of $m$ starting with $m=0$. As Table~\ref{tbl:candidate} shows, the value of
$E(181/437)$ is numerically smaller than $\lambda_{\max}(C_{s_*,283})$, so we
found a counterexample candidate $r=181/437$. We do not address whether $437$ is
the smallest possible denominator of a rational counterexample, since a single
counterexample suffices to disprove the conjecture.

\begin{table}[ht]
\begin{center}
\begin{tabular}{r@{\quad}r@{\quad}c|r@{\quad}c}
\hline
$p$ & $q$ & $E(p/q)$ & radius $m$ & $\lambda_{\max}(C_{s_*,m})$ \\
\hline
$1$   & $3$   & $2.732050807569$ & $0$ & $2.000000000000$ \\
$2$   & $5$   & $2.618033988750$ & $1$ & $2.476942779515$ \\
$3$   & $7$   & $2.610973927190$ & $2$ & $2.549179455498$ \\
$5$   & $12$  & $2.594619588219$ & $3$ & $2.575314207478$ \\
$7$   & $17$  & $2.593010334354$ & $4$ & $2.579669562737$ \\
$12$  & $29$  & $2.591391623995$ & $5$ & $2.582621869070$ \\
$17$  & $41$  & $2.591356427345$ & $100$ & $2.591038401762$ \\
$19$  & $46$  & $2.591309822224$ & $120$ & $2.591043375837$ \\
$29$  & $70$  & $2.591101474341$ & $280$ & $2.591046746750$ \\
$41$  & $99$  & $2.591070020988$ & $282$ & $2.591046772718$ \\
$70$  & $169$ & $2.591051417792$ & \cellcolor{paperblue} $283$ &
\cellcolor{paperblue} $2.591046796662$ \\
$111$ & $268$ & $2.591048970718$ & \cellcolor{paperblue} $284$ &
\cellcolor{paperblue} $2.591046806337$ \\
$169$ & $408$ & $2.591048315013$ & \cellcolor{paperblue} $288$ &
\cellcolor{paperblue} $2.591046917842$ \\
\cellcolor{paperblue} $181$ & 
\cellcolor{paperblue}$437$ & 
\cellcolor{paperblue}$2.591046792069$ & \cellcolor{paperblue} $300$ &
\cellcolor{paperblue} $2.591047022536$ \\
\hline
\end{tabular}
\end{center}
\caption{The search for a candidate counterexample. All decimal numbers are approximate.}\label{tbl:candidate}
\end{table}

The displayed values in this table are numerical estimates, not rigorous bounds.
Their only purpose is to show how the candidate was discovered. The candidate is
strictly smaller than the silver mean. Indeed
\begin{equation}
   \frac{181}{437}+1=\frac{618}{437}
\end{equation}
and
\begin{equation}
   2\cdot 437^2-618^2=14>0.
\end{equation}
Hence
\begin{equation}
   \frac{181}{437}<\sqrt2-1.
\end{equation}
The numerical gap between $E(181/437)$ and $\lambda_{\max}(C_{s_*,283})$ is
approximately $4.6\times10^{-9}$. Once the candidate $r=181/437$, the finite
radius $m=283$, and a suitable threshold $2.591046794$ had been identified
numerically, we show below how to certify the inequality
\begin{equation}
   E(181/437)=\lambda_{\max}(H_{181/437}) < 2.591046794
   < \lambda_{\max}(C_{s_*,283})\leqslant E(s_*).
\end{equation}
The comparison is certified using only exact integer arithmetic, without relying
on the floating-point eigenvalue estimates that suggested it. This proves that
the silver mean conjecture is wrong by explicitly presenting a (rational)
$\alpha\equiv r$ whose Harper norm $E(r)$ is strictly smaller than $E(s_*)$.

\section{Two elementary positive-vector bounds}

Here we establish two simple results that allow us to certify the inequalities
for $E(\alpha)$ with integer arithmetic only. The first lemma provides a
certificate for the upper bound and is valid for rational $\alpha$. The second
lemma gives a certificate for the lower bound and is valid for all $\alpha$.

Let $A=(a_{ij})$ be a real symmetric $n\times n$ matrix such that
$a_{ij}\geqslant 0$ for $i\ne j$, and suppose that the graph with vertices
$1,\dots,n$ and edges $(i, j)$ with $a_{ij}>0$ is connected. The following
lemma shows how to verify the bound on the maximal eigenvalue without computing
it first. Inequalities between vectors in this section are componentwise.

\begin{lemma}\label{lem:upper} Let $x\in\mathbb R^n$ have non-negative
components. If, for some $t\in\mathbb R$
\begin{equation}\label{eq:positive-upper}
 Ax \leqslant t x,
\end{equation}
componentwise, and the inequality is strict in at least one coordinate, then
\begin{equation}\label{eq:At}
    \lammax(A) < t.
\end{equation}
\end{lemma}
\begin{proof}
Because $A$ has nonnegative off-diagonal entries and its associated graph is
connected, the top eigenvalue of $A$ has a strictly positive eigenvector.
Choose $c>0$ sufficiently large that $B=A+cI$ is entrywise nonnegative, where
$I$ is the identity matrix. The matrix $B$ is irreducible because its
off-diagonal graph is connected. Hence Perron-Frobenius gives an eigenvector
$v>0$ corresponding to the largest eigenvalue of $B$. Since $B=A+cI$, the
same vector satisfies 
\begin{equation}
    Av = \lammax(A)v.
\end{equation}
Now, by assumption $Ax \leqslant tx$ componentwise, with strict inequality in at
least one coordinate. Since every coordinate of $v$ is strictly positive,
multiplication by $v_j$ and summation gives the strict inequality
\begin{equation}
    v^TAx<t\,v^Tx.
\end{equation}
Because $A$ is symmetric, $v^TAx=(Av)^Tx$. Using
\begin{equation}
Av=\lambda_{\max}(A)v,
\end{equation}
we obtain
\begin{equation}
\lambda_{\max}(A)\,v^Tx<t\,v^Tx.
\end{equation}
Finally, since the inequality \eqref{eq:positive-upper} is strict in at least
one coordinate, $x\neq0$. As $x\geqslant 0$ and $v>0$, we have $v^Tx>0$.
Dividing by $v^Tx$ yields the inequality \eqref{eq:At}.
\end{proof}

\begin{lemma}\label{lem:lower} Let $\alpha,t\in\mathbb R$, let $m\geqslant1$ be an
integer, and let $y=(y_j)$ be a real sequence such that $y_j>0$
for $-m\leqslant j\leqslant m$, with $y_j=0$ outside this interval. Suppose that
\begin{equation}\label{eq:finite-subsolution}
 y_{j-1}+y_{j+1}+2\cos(2\pi\alpha j)y_j\geqslant t y_j
\end{equation}
for $-m\leqslant j\leqslant m$ and at least one of these inequalities is strict.
Then $E(\alpha)>t$. It follows that any vector $y$ that satisfies these
conditions certifies the specified lower bound.
\end{lemma}
\begin{proof}
Multiply each inequality by the positive number $y_j$ and sum over the
support. Since $y$ vanishes outside the interval,
\begin{equation}
 \ip{y}{H_{\alpha,0}y}
 =\sum_{j=-m}^{m}y_j(H_{\alpha,0}y)_j
 >t\sum_{j=-m}^{m}y_j^2=t\norm{y}^2.
\end{equation}
Thus $E(\alpha)\geqslant e(H_{\alpha,0})>t$. No inequalities outside the support
are needed, because their contributions to the inner product are multiplied by
zero.
\end{proof}

\section{Bounding the coefficients}\label{sec:recurrences}

We will prove an inequality of the form $E(r)<t<E(s)$, where $r,t$ are rational
and $s$ is real, by applying the two lemmas above. In the concrete application,
$s\equiv s_*=\sqrt2-1$ is irrational. In general, the values
$c_j(\alpha)=2\cos(2\pi\alpha j)$ are not rationals, and it is this fact that
requires an additional construction before we can apply the lemmas. For a
rational $\alpha$ these values are algebraic, but even in this case precise
comparisons are problematic if the only tool we use is elementary arithmetic.
The idea is to bound $c_j(\alpha)$ by rational numbers and use these bounds to
find the certificate sequences $x$ and $y$ from the lemmas.

Let $t = T/D$, where $T$ and $D>0$ are integers to be determined. We bound
$c_j(r)$ and $c_k(s)$ by rational numbers with the same denominator $D$, so we
are looking for the inequalities of the form
\begin{equation}\label{eq:ul}
 c_j(r)<\frac{u_j}{D}, \quad \frac{\ell_k}{D}<c_k(s).
\end{equation}
We need only a finite number of these bounds. For a reduced rational number
$r=p/q$ with $q \geqslant 3$, $\gcd(p,q)=1$, define $h=\floor{q/2}$. In
\eqref{eq:ul}, the index $j$ ranges over $0\leqslant j\leqslant h$, while $k$
ranges over $-m\leqslant k\leqslant m$, where $m\geqslant1$ is an integer chosen
independently of $r$. By evenness of the cosine, we may set $\ell_{-k}=\ell_k$.
We show how to choose $T$, $D$, $m$ and construct these numerators later.

Assuming that these parameters have already been found, we now construct the
main ingredient of the two witnesses $x$ and $y$ from the lemmas. Starting with
$a_0=D$, define
\begin{equation}\label{eq:upper-start}
 a_1=\floor{\frac{T-u_0}{2}},
\end{equation}
and then
\begin{equation}\label{eq:upper-recurrence}
 a_{j+1}=\floor{\frac{(T-u_j)a_j}{D}}-a_{j-1},
 \qquad 1\leqslant j\leqslant h.
\end{equation}
This produces $a_0,\ldots,a_h$ and one auxiliary terminal value $a_{h+1}$. The
witness itself will use only $a_0,\ldots,a_h$. The direction of rounding
enforces the required local inequalities. Indeed, $a_0=D$ and
\eqref{eq:upper-start} give
\begin{equation}\label{eq:upper-center}
 2a_1+\frac{u_0}{D}a_0\leqslant t a_0.
\end{equation}
For $1\leqslant j\leqslant h$, Eq.~\eqref{eq:upper-recurrence} gives
\begin{equation}\label{eq:upper-local}
 a_{j-1}+a_{j+1}+\frac{u_j}{D}a_j\leqslant t a_j.
\end{equation}
This is an exact inequality, not a statement about the numerical stability
of a recurrence.

Starting with $b_0=D$, define
\begin{equation}\label{eq:lower-start}
 b_1=\ceil{\frac{T-\ell_0}{2}},
\end{equation}
and then
\begin{equation}\label{eq:lower-recurrence}
 b_{j+1}=\ceil{\frac{(T-\ell_j)b_j}{D}}-b_{j-1},
 \qquad 1\leqslant j\leqslant m.
\end{equation}
This produces $b_0,\ldots,b_m$ and one auxiliary terminal value $b_{m+1}$.
Rounding up gives
\begin{equation}
 2b_1+\frac{\ell_0}{D}b_0\geqslant t b_0, \quad
 b_{j-1}+b_{j+1}+\frac{\ell_j}{D}b_j\geqslant t b_j, \quad
 1\leqslant j\leqslant m.\label{eq:lower-local}
\end{equation}
The auxiliary values $a_{h+1}$ and $b_{m+1}$ are sometimes called
\emph{ghost values}: they are used to check the eventual boundary
inequality, but are not entries of the final witnesses. Computing them
uses only the enclosures at sites $h$ and $m$; no potential at site
$h+1$ or $m+1$ is required.

\section{The general positive-sequence certificate}\label{sec:closure}

\begin{proposition}\label{prop:certificate} Let $r=p/q$ be a reduced rational
with $q\geqslant3$, put $h=\lfloor q/2\rfloor$, and let $s\in\mathbb R$. Let
$m,D,T$ be integers with $m\geqslant1$ and $D>0$, and set $t=T/D$. Suppose the
integer numerators $u_0,\ldots,u_h$ and $\ell_0,\ldots,\ell_m$ satisfy the
strict enclosures \eqref{eq:ul}. Generate $a_0,\ldots,a_{h+1}$ and
$b_0,\ldots,b_{m+1}$, starting from $a_0=b_0=D$, by
\eqref{eq:upper-start}--\eqref{eq:upper-recurrence} and
\eqref{eq:lower-start}--\eqref{eq:lower-recurrence}. Assume the upper sign
conditions
\begin{equation}\label{eq:upper-signs}
 a_j>0\quad(0\leqslant j\leqslant h),\qquad
 \begin{cases}
  a_{h+1}\geqslant a_h,&q=2h+1,\\
  a_{h+1}\geqslant a_{h-1},&q=2h,
 \end{cases}
\end{equation}
and the lower sign conditions
\begin{equation}\label{eq:lower-signs}
 b_j>0\quad(0\leqslant j\leqslant m),\quad b_{m+1}\leqslant 0.
\end{equation}
Then
\begin{equation}\label{eq:general-conclusion}
    E(r)<t<E(s).
\end{equation}
\end{proposition}

\begin{proof}
We construct two witnesses $x$ and $y$ from the sequences $a$ and $b$. The
vector $x$ has $q$ elements, while the sequence $y$ is supported on $2m+1$
sites. Define a $q$-vector $x$ by
\begin{equation}\label{eq:periodic-witness}
 x_j=a_{\min(j,q-j)},\qquad 0\leqslant j<q.
\end{equation}
All its coordinates are positive. The coefficients $c_j(r)$ are symmetric:
\begin{equation}
 c_{q-j}(r)=2\cos(2\pi p-2\pi pj/q)=c_j(r).
\end{equation}
At site zero, both neighbors equal $a_1$, so
\eqref{eq:upper-center} applies. At the sites corresponding to
$1\le j<h$, the neighboring values are $a_{j-1}$ and $a_{j+1}$, so
\eqref{eq:upper-local} applies. It remains to verify the terminal sites.

\medskip\noindent Odd period: $q=2h+1$. The vector has the form
\begin{equation}
 x = (a_0,a_1,\ldots,a_{h-1},a_h,a_h,a_{h-1},\ldots,a_1).
\end{equation}
For $q=3$ this simply means $(a_0,a_1,a_1)$. At either of the two terminal
sites with value $a_h$, the actual neighbors have values $a_{h-1}$ and
$a_h$. But the recurrence at $j=h$ used $a_{h+1}\geqslant a_h$.
Therefore
\begin{equation}\label{eq:closing}
 a_{h-1}+a_h+\frac{u_h}{D}a_h
 \leqslant a_{h-1}+a_{h+1}+\frac{u_h}{D}a_h
 \leqslant t a_h.
\end{equation}
This verifies the bond between the two terminal sites.

\medskip\noindent Even period: $q=2h$. Here $h\geqslant2$, and the vector has
the form
\begin{equation}
 x = (a_0,a_1,\ldots,a_{h-1},a_h,a_{h-1},\ldots,a_1).
\end{equation}
There is a single terminal site with value $a_h$, and both its neighbors
have value $a_{h-1}$. Thus the condition needed to close the cycle is
$a_{h+1}\geqslant a_{h-1}$, rather than $a_{h+1}\geqslant a_h$. It gives
\begin{equation}\label{eq:closing-even}
 2a_{h-1}+\frac{u_h}{D}a_h
 \leqslant a_{h-1}+a_{h+1}+\frac{u_h}{D}a_h
 \leqslant t a_h.
\end{equation}
In particular, for $q=4$ the witness is $(a_0,a_1,a_2,a_1)$, and the terminal
check is $a_3\geqslant a_1$.

\medskip\noindent In either parity case, the enclosed periodic operator applied
to $x$ is at most $t x$. Replacing $u_j/D$ by the strictly smaller true
coefficient $c_j(r)$ decreases every component strictly, because each $x_j$ is
positive. Consequently,
\begin{equation}
    H_{p/q}x<tx
\end{equation}
at all $q$ sites. Applying Lemma~\ref{lem:upper} with $A=H_{p/q}$ and then using
\eqref{eq:rational}, we obtain
\begin{equation}
    E(r)=\lambda_{\max}(H_{p/q})<t.
\end{equation}

Now define a finitely supported sequence on $\Z$ by
\begin{equation}\label{eq:dirichlet-witness}
 y_k=\begin{cases}
 b_{|k|},&|k|\leqslant m,\\
 0,&|k|>m.
 \end{cases}
\end{equation}
It is positive on its $2m+1$-site support. The cosine is even, so it suffices to
check nonnegative indices. At zero and at $1\leqslant k<m$, use
\eqref{eq:lower-local}. The interior range is empty when $m=1$, and the center
argument is unchanged. At the right endpoint, the recurrence gives
\begin{equation}
 b_{m-1}+b_{m+1}+\frac{\ell_m}{D}b_m\geqslant t b_m.
\end{equation}
Since $b_{m+1}\leqslant 0$, replacing this auxiliary value by the actual
exterior value $0$ does not decrease the left side:
\begin{equation}\label{eq:cutoff}
 b_{m-1}+\frac{\ell_m}{D}b_m
 \geqslant t b_m-b_{m+1}\geqslant t b_m.
\end{equation}
The same holds at the left endpoint. Replacing each $\ell_k/D$ by the larger
true coefficients $c_k(s)$ makes the required inequalities strict. Hence
\begin{equation}
 (H_{s,0}y)_k > t y_k
\end{equation}
for $-m\leqslant k\leqslant m$. Lemma~\ref{lem:lower} gives $E(s)>t$, completing
the proof.
\end{proof}

The proposition is a sufficient-condition statement. It does not assert that the
checks succeed for every choice of parameters, and failure of a sign check does
not establish the opposite spectral ordering. An arbitrary real $\alpha\equiv s$
is allowed, but an arithmetic application requires proved enclosures for its
finitely many coefficient values $c_k(s)$. 

\section{Specialization to the silver-mean counterexample}\label{sec:application}

We now fix the parameters. Throughout this and the next sections, fix
\begin{equation}\label{eq:params}
    p=181,\qquad q=437,\qquad r=\frac{181}{437},
    \qquad s=s_*=\sqrt2-1,
    \qquad h=218,\qquad m=283.
\end{equation}
The preceding proposition and its proof did not use
them. We use the radius $m=283$ that already revealed the candidate in
Table~\ref{tbl:candidate}. The separating threshold $t$ is chosen inside the
narrow gap suggested by that table. The numerical estimates motivate this
choice, but none of them is used in the exact verification below. Take
\begin{equation}\label{eq:parameters}
 D=10^{10},\quad T=25910467940,\quad t=\frac TD=2.591046794.
\end{equation}
Since $q=437$ is odd, the sign conditions required by
Proposition~\ref{prop:certificate} become
\begin{equation}\label{eq:concrete-signs}
 a_0,\ldots,a_{218}>0,\quad a_{219}\geqslant a_{218}, \quad
 b_0,\ldots,b_{283}>0,\quad b_{284}\leqslant0.
\end{equation}

\begin{theorem}\label{thm:main} For the parameters \eqref{eq:parameters} the
inequalities \eqref{eq:intro-counterexample} are valid. In particular, the
silver mean is not a global minimizer of the Harper norm.
\end{theorem}

The periodic upper witness is explicitly
\begin{equation}
 x=(a_0,a_1,\ldots,a_{217},a_{218},a_{218},a_{217},\ldots,a_1),
\end{equation}
with $437$ entries. The lower witness is $y_j=b_{|j|}$ for $|j|\leqslant283$ and
zero otherwise, so it has $567$ nonzero entries. Their auxiliary values
$a_{219}$ and $b_{284}$ are used only for the terminal tests. It remains to
construct the integer numerators in \eqref{eq:ul} and verify
\eqref{eq:concrete-signs}. Below we prove the rational Taylor bounds used to
generate those numerators; Section~\ref{sec:arithmetic} gives the exact sequence
values and sign checks that verify the hypotheses of
Proposition~\ref{prop:certificate}. All subsequent uses of $p,q,r,s,h,m,D,T$
refer to this specialization. For these concrete parameters, the remaining
analytic task is to construct the integers $u_j$ and $\ell_j$. We use a single
denominator
\begin{equation}\label{eq:Q}
 Q=10^{16}
\end{equation}
for approximating both $\pi$ and $s_*$. All polynomial evaluations below are
performed as exact rational arithmetic.

\subsection{A rational approximation to \texorpdfstring{$\pi$}{pi}}

Let
\begin{equation}\label{eq:atan}
 S(q)=\sum_{k=0}^{19}\frac{(-1)^k}{(2k+1)q^{2k+1}}.
\end{equation}
The alternating-series estimate gives
\begin{equation}
 \left|\arctan(1/q)-S(q)\right|
 <\frac{1}{41q^{41}}\qquad(q>1).
\end{equation}
Machin's identity \cite{WeissteinMachin} is
\begin{equation}\label{eq:machin}
 \pi=16\arctan(1/5)-4\arctan(1/239).
\end{equation}
Set
\begin{equation}\label{eq:pi0}
 \widehat\pi=16S(5)-4S(239),\qquad
 \pi_0=\frac{\floor{Q\widehat\pi}}{Q}.
\end{equation}
The exact rational comparison
\begin{equation}\label{eq:pi-remainder}
 \frac{16}{41\cdot5^{41}}+\frac4{41\cdot239^{41}}<\frac1Q
\end{equation}
implies $|\widehat\pi-\pi|<1/Q$. Flooring adds an error of less than
$1/Q$, hence
\begin{equation}\label{eq:pi-bound}
 |\pi-\pi_0|<\frac2Q.
\end{equation}
We have
\begin{equation}
 \pi_0=\frac{31415926535897932}{10^{16}}.
\end{equation}
This displayed rational is generated by \eqref{eq:atan}--\eqref{eq:pi0};
it is not an unproved approximation.

\subsection{A rational approximation to the silver mean}

For a nonnegative integer $N$, write
\begin{equation}
    \operatorname{isqrt}(N)=\lfloor\sqrt N\rfloor
\end{equation}
for its integer square root. Set
\begin{equation}
    k=\operatorname{isqrt}(2Q^2),\qquad
    s_0=\frac{k-Q}{Q}.
\end{equation}
The computed integer $k$ is verified by checking
\begin{equation}\label{eq:isqrt}
    k\geqslant0,\qquad k^2\leqslant2Q^2<(k+1)^2.
\end{equation}
Dividing by $Q^2$ and taking nonnegative square roots gives
\begin{equation}\label{eq:s-bound}
 0\leqslant s_*-s_0<\frac1Q.
\end{equation}
Here the computed values are
\begin{equation}
 k=14142135623730950,\quad
 s_0=\frac{4142135623730950}{10^{16}}.
\end{equation}
It is the condition \eqref{eq:isqrt} that establishes the inequality
\eqref{eq:s-bound}, rather than a floating-point approximation to the square
root.

\subsection{Taylor evaluation after exact reduction modulo one}

For the rational sequence, set $(\alpha,\alpha_0)=(r,r)$; for the irrational
sequence, set $(\alpha,\alpha_0)=(s_*,s_0)$. At index $j$, define
\begin{equation}\label{eq:angle}
 n_j=\floor{j\alpha_0+\frac12},\qquad
 v_j=j\alpha_0-n_j\in[-1/2,1/2),\qquad \xi_j=2\pi_0v_j.
\end{equation}
These quantities are exact rationals or integers. Since $0<\pi_0<4$,
we have $|\xi_j|<4$.

Define the degree-$24$ polynomial
\begin{equation}\label{eq:taylor}
 B(x)=2\sum_{k=0}^{12}\frac{(-1)^kx^{2k}}{(2k)!}.
\end{equation}
This is also the Taylor polynomial through degree $25$, because the
coefficient of $x^{25}$ is zero. Taylor's theorem, with the $26$th
derivative bounded in absolute value by $2$, therefore gives
\begin{equation}\label{eq:taylor-bound}
 |2\cos x-B(x)|\leqslant\frac{2|x|^{26}}{26!}
 <\frac{2\cdot4^{26}}{26!}
 \qquad(|x|<4).
\end{equation}
To compare the exact target angle with $\xi_j$, subtract the same integer
$n_j$ from $j\alpha$. This avoids any assumption that nearby frequencies
fall in the same nearest-integer cell. Periodicity gives
\begin{equation}
 c_j(\alpha)=2\cos\bigl(2\pi(j\alpha-n_j)\bigr),
\end{equation}
and
\begin{equation}
 |2\pi(j\alpha-n_j)-2\pi_0v_j|
 \leqslant2\pi j|\alpha-\alpha_0|+2|\pi-\pi_0||v_j|
 <\frac{8j+2}{Q},
\end{equation}
using $\pi<4$, \eqref{eq:pi-bound} and \eqref{eq:s-bound}. The same bound holds
on the rational side, where $\alpha-\alpha_0=0$. Since the derivative of $2\cos
x$ has absolute value at most $2$, it follows that
\begin{equation}\label{eq:angle-error}
 |c_j(\alpha)-2\cos \xi_j|<\frac{16j+4}{Q}.
\end{equation}
The rational upper construction requires coefficient enclosures
for $0\leqslant j\leqslant h=218$, while the lower construction
requires them for $0\leqslant j\leqslant m=283$.
We therefore use a common uniform estimate on
\begin{equation}
    0\leqslant j\leqslant\max(h,m)=283.
\end{equation}
The auxiliary values $a_{h+1}$ and $b_{m+1}$ require only
the coefficients at indices $h$ and $m$, respectively,
so no larger coefficient index is needed.
Combining \eqref{eq:taylor-bound} and \eqref{eq:angle-error}, uniformly for
$0\leqslant j\leqslant283$, yields
\begin{equation}\label{eq:eps}
 |c_j(\alpha)-B(\xi_j)|<\eps,
 \qquad
 \eps=\frac{16\cdot283+4}{Q}+\frac{2\cdot4^{26}}{26!}<\frac1D.
\end{equation}
The last comparison is checked exactly. For example, in reduced rational
form it is
\begin{equation}
 D\eps=
 \frac{175284785173865693551}{769217417004786750000}<1.
\end{equation}

\subsection{Outward rounding to integer potential numerators}

Put
\begin{equation}\label{eq:z}
 z_j^{(\alpha)}=\floor{D B(\xi_j)}.
\end{equation}
The definition of the floor and \eqref{eq:eps} imply
\begin{equation}\label{eq:cos-enclosure}
 \frac{z_j^{(\alpha)}-1}{D}
 <c_j(\alpha)
 <\frac{z_j^{(\alpha)}+2}{D}.
\end{equation}
Indeed, $z_j^{(\alpha)}/D\leqslant B(\xi_j)<(z_j^{(\alpha)}+1)/D$, and the
additional error is strictly smaller than $1/D$ on either side. Thus the
required definitions are
\begin{equation}\label{eq:integer-potentials}
 u_j=z_j^{(r)}+2\quad(0\leqslant j\leqslant218),\qquad
 \ell_j=z_j^{(s_*)}-1\quad(0\leqslant j\leqslant283).
\end{equation}

At $j=0$ this uniform rule gives
\begin{equation}\label{eq:zero-potentials}
 u_0=20000000002=2D+2,\qquad
 \ell_0=19999999999=2D-1.
\end{equation}
Although the true value of $c_0(\alpha)$ is exactly $2$, these harmless outward
offsets are retained so that one formula covers every site. The proof of
strictness on the rational side uses the strict enclosures in
\eqref{eq:cos-enclosure}.

\section{The exact arithmetic certificate}\label{sec:arithmetic}

The complete scalar specification is given by equations \eqref{eq:params} and
\eqref{eq:Q} together with $20$ terms in each arctangent sum and the $13$ even
terms in \eqref{eq:taylor}. There are no precomputed cosine arrays or witness
vectors to import. Both sequences $a$ and $b$ are generated using the threshold
numerator $T=25910467940$ in
\eqref{eq:upper-start}--\eqref{eq:lower-recurrence}. Running the constructions
gives the following exact integers:
\begin{equation}\label{eq:upper-result}
\begin{split}
 \min_{0\leqslant j\leqslant218}a_j&=61141958>0, \quad
 a_{218}=61141958,\quad a_{219}=61145499>a_{218}, \\
 \min_{0\leqslant j\leqslant283}b_j&=1727214>0, \quad
 b_{283}=1727214,\quad b_{284}=-151774<0.
\end{split}
\end{equation}
These minima summarize checks of every entry of the relevant sequences.

At site zero, \eqref{eq:zero-potentials} gives
\begin{equation}
 a_1=\floor{\frac{5910467938}{2}}=2955233969,\qquad
 b_1=\ceil{\frac{5910467941}{2}}=2955233971.
\end{equation}
The different directions of rounding already appear here. At the last upper step
the computed data are
\begin{equation}
 u_{218}=-5326676596,\qquad
 a_{217}=129844518,\qquad a_{218}=61141958.
\end{equation}
Substituting into \eqref{eq:upper-recurrence} gives
\begin{equation}
 a_{219}=\floor{\frac{31237144536\cdot61141958}{10^{10}}}
       -129844518 =190990017-129844518=61145499.
\end{equation}
In particular, the cyclic boundary can be closed because this value
exceeds $a_{218}$.

At the last lower step the data are
\begin{equation}
 \ell_{283}=3446238113,\qquad
 b_{282}=4031828,\qquad b_{283}=1727214.
\end{equation}
Therefore
\begin{equation}
 b_{284}=\ceil{\frac{22464229827\cdot1727214}{10^{10}}}-4031828
 =3880054-4031828=-151774.
\end{equation}
Every floor or ceiling can be checked by multiplication with its positive
denominator. All these computations are performed in the supplementary
Mathematica file \texttt{silver\_mean\_certificate.wl}. The coefficient enclosures and
the sign checks above satisfy the hypotheses of
Proposition~\ref{prop:certificate}, which proves Theorem~\ref{thm:main}.

\appendix

\section{Proof of Proposition~\ref{prop:dirichlet-convergence}}

As observed above, a top eigenvector of $B_{\alpha,m}$ can be chosen
even. Hence the largest eigenvalue of $B_{\alpha,m}$ is already
attained on the even subspace, and therefore
\begin{equation}
    \lambda_{\max}(C_{\alpha,m})
    =
    \lambda_{\max}(B_{\alpha,m}).
\end{equation}
If $x$ is a real vector supported on $[-m,m]$, extend it by zero to a vector
$\widetilde x\in\ell^2(\mathbb Z)$. Real vectors suffice because the operator
has real coefficients. Then $\|\widetilde x\|=\|x\|$. Every hopping term crossing
the boundary contains an exterior zero coordinate, so
\begin{equation}
    \langle \widetilde x,H_{\alpha,0}\widetilde x\rangle
    =
    x^TB_{\alpha,m}x.
\end{equation}
Consequently
\begin{equation}
    \lambda_{\max}(B_{\alpha,m})
    \leqslant
    \sup_{0\ne z\in\ell^2(\mathbb Z)}
    \frac{\langle z,H_{\alpha,0}z\rangle}{\|z\|^2}.
\end{equation}
By \eqref{eq:intro-zero-phase}, the right-hand side is
$e(H_{\alpha,0})=E(\alpha)$. Hence
\begin{equation}
    \lambda_{\max}(C_{\alpha,m})\leqslant E(\alpha).
\end{equation}
It remains to prove convergence. The intervals $[-m,m]$ are nested. Hence every
vector supported on $[-m,m]$ is also a valid trial vector for $B_{\alpha,m+1}$,
and therefore
\begin{equation}
    \lambda_{\max}(B_{\alpha,m})
    \leqslant
    \lambda_{\max}(B_{\alpha,m+1}).
\end{equation}
Thus the sequence $\lambda_{\max}(C_{\alpha,m})$ is nondecreasing and bounded
above by $E(\alpha)$.

Now let $\varepsilon>0$. By the variational characterization of the
upper spectral edge, there exists a nonzero
$z\in\ell^2(\mathbb Z)$ such that
\begin{equation}
    \frac{\langle z,H_{\alpha,0}z\rangle}{\|z\|^2}
    >
    E(\alpha)-\frac{\varepsilon}{2}.
\end{equation}
Finitely supported vectors are dense in $\ell^2(\mathbb Z)$, and
$H_{\alpha,0}$ is bounded. Hence we may choose a finitely supported
vector $y$ such that
\begin{equation}
    \frac{\langle y,H_{\alpha,0}y\rangle}{\|y\|^2}
    >
    E(\alpha)-\varepsilon.
\end{equation}
Choose $m$ large enough that $\operatorname{supp}y\subset[-m,m]$.
Since extension by zero preserves the quadratic form,
\begin{equation}
    \frac{y^TB_{\alpha,m}y}{\|y\|^2}
    =
    \frac{\langle y,H_{\alpha,0}y\rangle}{\|y\|^2}
    >
    E(\alpha)-\varepsilon.
\end{equation}
Therefore
\begin{equation}
    \lambda_{\max}(C_{\alpha,m})
    =
    \lambda_{\max}(B_{\alpha,m})
    >
    E(\alpha)-\varepsilon.
\end{equation}
Together with
\begin{equation}
    \lambda_{\max}(C_{\alpha,m})\leqslant E(\alpha),
\end{equation}
this proves
\begin{equation}
    \lambda_{\max}(C_{\alpha,m})
    \longrightarrow E(\alpha).
\end{equation}
Since the sequence is nondecreasing, we may write
\begin{equation}
    \lambda_{\max}(C_{\alpha,m})
    \nearrow E(\alpha),
\end{equation}
which concludes the proof. \\

\noindent The authors thank the BMFTR (via QR.N, QuaPhySI, and QuKuK) and the EU
(via CLUSTEC) for their support.

\end{document}